\def\sqr#1#2{{\vcenter{\hrule height.#2pt
      \hbox{\vrule width.#2pt height#1pt \kern#1pt
         \vrule width.#2pt}
       \hrule height.#2pt}}}

\def\th #1 #2. #3\par{\medbreak{\bf#1 #2.
\enspace}{\sl#3}\par\medbreak}
\documentclass{llncs1}
\usepackage{setspace}
\usepackage{makeidx}  
\usepackage{graphicx}
\begin{document}
\frontmatter          
\pagestyle{headings}  

\title{On the zeros of the Riemann zeta function: Lessons to Learn from A False Proof of the Riemann Hypothesis}
\titlerunning{Lessons to Learn}  
%
\author{Jorma Jormakka}
\authorrunning{Jorma Jormakka}   
%
\tocauthor{Jorma Jormakka}

\institute{ Contact by: \email{jorma.o.jormakka@gmail.com}}

\maketitle              

\begin{abstract}
In 2008 I thought I found a proof of the Riemann Hypothesis, but there was an 
error. In the Spring 2020 I believed to have fixed the error, but it cannot 
be fixed. I describe here where the error was. It took me several 
days to find the error in a careful checking before a possible submission 
to a payable review offered by one leading journal. There were three simple 
lemmas and one simple theorem, all were correct, yet there was an error: 
what Lemma 2 proved was not exactly what Lemma 3 needed. So, it was the
connection of the lemmas. This paper came out empty, but I have found a 
different proof of the Riemann Hypothesis and it seems so far correct. 
In the discussion at the end of this paper
I raise a matter that I think is of importance to the review process in 
mathematics. 
\end{abstract}

Let $p_j$, $j=1,2,\dots$, denote the primes in the increasing order. 
The Riemann zeta
function $\zeta(s)$ (see e.g. [1] for basic facts of zeta)
satisfies for $Re\{s\}>1$ the equation
$$\zeta(s)^{-1}\zeta'(s)={d\over ds}\ln \zeta(s)=
-{d\over ds}\sum_{j=1}^{\infty}\ln(1-p_j^{-s})
=-\sum_{j=1}^{\infty}\ln(p_j)p_j^{-s}(1-p_j^{-s})^{-1}$$
$$=-\sum_{j=1}^{\infty}\ln(p_j)p_j^{-s}\left(1+\sum_{k=1}^{\infty}p_j^{-ks}\right)=h(s)+g(s)$$
where 
$$h(s)=-\sum_{j=1}^{\infty}\ln(p_j)p_j^{-s}$$
$$g(s)=-\sum_{j=1}^{\infty}\ln(p_j)p_j^{-2s}\sum_{k=0}^{\infty}p_j^{-ks}.$$
For any small $\rho>0$ and sufficiently large $j$
$$|p_j^{\rho}|>|\ln(p_j)\sum_{k=0}^{\infty}p_j^{-ks}|.$$
Therefore $g(s)$ converges absolutely if $Re\{s\}>{1\over 2}$. 
Let $s_0$ be a zero of zeta and ${1\over 2}<Re\{s_0\}<1$.
The series expression of $h(s)$ does not converge if 
${1\over 2}<Re\{s\}<1$, but $h(s)$ can be analytically continued to 
$Re\{s\}>{1\over 2}$ to all points where $\zeta(s)$ is not zero or infinte 
($s=1$) by 
$$h(s)=\zeta(s)^{-1}\zeta'(s)-g(s).$$
Thus, $h(s)$ has a finite value at all points $s$, $1>Re\{s\}>{1\over 2}$, 
where $\zeta(s)\not=0$.
From the equation
$$\zeta'(s)=h(s)\zeta(s)+g(s)\zeta(s)$$
we conclude that in zeros and the pole of $\zeta(s)$ with $Re\{s\}>{1\over 2}$
the function $h(s)$ has a simple pole and is of the type
$$h(s)={r\over s-s_0}+w(s)$$
where $s_0$ is the pole and $w(s)$ is analytic close to $s_0$. 
If it is possible to show that $h(s)$ is finite at $0<Re\{s\}<{1\over 2}$, 
then the Riemann Hypothesis is proven. What I tried is in a simplified
form as follows.

Let us define a function of two complex variables
$$\varphi(s,z)=\prod_{j=1}^\infty\left(1-z\ln(p_j)p_j^{-s-f(z)}\right).$$
Here $f(z)$ is some analytic function.
$\varphi(s,z)$ is analytic as a function of $s$ and as a function of $z$
in the area where the infinite product converges absolutely.
The function 
$\varphi(s,z)$ satisfies the equation if $Re\{s+f(z)\}>1$
$$\varphi(s,z)^{-1}{\partial\over \partial z}\varphi(s,z)
={\partial \over \partial z}\ln \varphi(s)=
{\partial\over \partial z}\prod_{j=1}^\infty\ln(1-z\ln(p_j)p_j^{-s-f(z)})$$
$$=\sum_{j=1}^n(-\ln(p_j)p_j^{-s-f(z)}+zf'(z)\ln(p_j)^2p_j^{-s-f(z)})
(1-z\ln(p_j)p_j^{-s-f(z)})^{-1}$$
$$=h(s+f(z))+zf'(z)h'(s+f(z))+u(s,z)$$
where 
$$u(s,z)=\sum_{j=1}^\infty p^{-2s-2f(z)}(-\ln(p_j)p_j^{-f(z)}+zf'(z)\ln(pj)^2)\sum_{k=0}^{\infty}(p_j^{-s-f(z)})^k.$$
For any small $\rho>0$ and sufficiently large $j$
$$|p_j^{\rho}|>|
(-\ln(p_j)p_j^{-f(z)}+zf'(z)\ln(pj)^2)\sum_{k=0}^{\infty}(p_j^{-s-f(z)})^k|.$$
Therefore $u(s,z)$ converges absolutely if $Re\{s+f(s)\}>{1\over 2}$.
As $u(s,z)$ converges and $h(s)$ was already analytically continued, 
the equation
$${\partial \over \partial z}\ln \varphi(s,z)=h(s+f(z))+zf'(z)h'(s+f(z))+u(s,z)$$
holds if $Re\{s+f(s)\}>{1\over 2}$.
We can integrate this equation
$$\ln \varphi(s,z)=zh(s+f(z))+\int^zu(s,z)dz+\phi(s).$$
Here $\phi(s)$ is the integration constant, an analytic function that does 
not depend on $z$. 
From the Taylor series of $\varphi(s,z)$ in $Re\{s+f(s)\}>1$.
$$\ln\varphi(s,z)=\sum_{j=1}^\infty\ln(1-z\ln(p_j)p_j^{-s-f(z)})$$
$$=z\sum_{j=1}^\infty\ln(p_j)p_j^{s+f(z)}+O(z^2p_j^{-2s-2f(z)})$$
we conclude that in an open environment in $Re\{s+f(s)\}>1$ there is
no term that does not depend on $z$. Therefore $\phi(s)$ is zero everywhere.
Thus,
$$\ln \varphi(s,z)=zh(s+f(z))+\int^zu(s,z)dz.$$
Close to the pole $s_0$ of $h(s)$ we get
$$\ln \varphi(s,z)=z{r\over s+f(z)-s_0}+w_1(s,z)$$
where $w_1(s,z)$ is finite. Consequenty 
$$\ln \varphi(s_0,z)=r{z\over f(z)}+w_1(s_0,z).$$
Depending on how we select $f(s)$ we can have $\ln\varphi(s,z)$ as a finite
function with zero value, or nonzero finite value, or infinite value. 

Let $s=s_0-z$. Then for 
$$\varphi(s,z)=\prod_{j=1}^\infty\left(1-z\ln(p_j)p_j^{-s}\right)$$
holds
$${\partial \over \partial z}\ln \varphi(s,z)=h(s)+u(s,z)\eqno(1)$$ 
while for
$$\varphi_1(s_0,z)=\prod_{j=1}^\infty\left(1-z\ln(p_j)p_j^{-s_0+z}\right)$$
holds
$${\partial \over \partial z}\ln \varphi_1(s_0,z)=h(s_0-z)-zh'(s_0-z)+u_1(s_0,z).$$
Clearly
$$\ln \varphi(s,z)=\ln \varphi_1(s_0,z).$$
Actually I used in 2008 a slightly different function. There I defined
$$\varphi_1(s,z)=\prod_{j=1}^\infty (1-p_j^{-s}+p_j^{-s-z}).$$
The difference is not essential since
$$1-p_j^{-s}+p_j^{-s-z}=1-z\ln(p_j)p_j^{-s}+O(z^2\sum_{j=1}^\infty p_j^{-s}).$$
The term $\sum_{j=1}^\infty p_j^{-s}$ in the $O()$ term 
has a simple pole at $s_0$, but 
$z^2$ takes it to zero when $z\to 0$. With this $\varphi_1(s,z)$ 
the equation (1) takes the form 
$${\partial \over \partial z}\ln \varphi_1(s,z)=h(s+z)+u_2(s,z)$$
where $u_2(s,z)$ is finite. The equation can be written as 
$${\partial \over \partial z}\varphi_1(s,z)=h(s+z)\varphi_1(s,z)+u_2(s,z)\varphi_1(s,z).\eqno(3)$$
As $s=s_0-z$ the term $h(s+z)=h(s_0)$. It seemed to me that it is sufficient
to show that there is a point $(s,z)$, $s+z=s_0$, such that $\varphi_1(s,z)$
is finite and nonzero. Then the right side is infinite if $h(s_0)$ is 
infinite. It happens to be that $\varphi_1(s,z)$ is finite and nonzero.
This is because 
$$\varphi_1(s,z)=\varphi_2(s_0,z)+O(z^2\sum_{j=1}^\infty p_j^{-s}).$$ 
The $O(z^2\sum_{j=1}^\infty p_j^{-s})$ term is finite even at $z=0$
and $\varphi_2(s_0,z)$ is finite for any small $z$, including $z=0$. 
The partial differential in the left side appeared to me obviously 
finite, as the function $\varphi_1(s,z)$ was finite and continuous. 
But it is here where the error was. Let us see this error.
For any function $g(s,z)=g(s',z')$ the calculation
$${\partial g\over \partial z}={ds'\over dz}{\partial g\over \partial s'}
+{dz'\over dz}{\partial g\over \partial z'}$$
is correct.
As an example, if $s=s'-z$, $z=z'$ and $g(s,z)=sz$, then $g(s',z')=(s'-z')z'$.
Thus
$${\partial g(s,z)\over \partial z}=s$$
and
$${\partial g(s',z')\over \partial z}
={ds'\over dz}{\partial g(s',z')\over \partial s'}+{dz'\over dz}{\partial g(s',z')\over \partial z'}$$
$$=1\cdot z'+1\cdot (s'-2z')=z'+s'-2z'=s'-z'=s'-z=s$$
is correct. Here holds
$${\partial g(s,z)\over \partial z}={\partial g(s',z')\over \partial z}.$$
Let us do our case in a similar way. We set 
$s'=s-z$ and $z'=z$ and define
$$g(s,z)=\ln(\prod_{j=1}^\infty (1-z\ln(p_j)p^{-s}))=\ln\varphi(s,z).$$ 
Then simply inserting $s$ and $z$ by $s=s'-z'$, $z'=z$ we get
$$g(s',z')=\ln(\prod_{j=1}^\infty (1-z'\ln(p_j)p^{-s'+z}))$$ 
$$=\ln(\prod_{j=1}^\infty (1-z\ln(p_j)p^{-s'+z}))=\ln \varphi_2(s',z).$$ 
We derived above
$${\partial g(s,z)\over \partial z}=h(s)+u(s,z)$$
and 
$${\partial g(s',z')\over \partial z'}=h(s'-z')-z'h'(s'-z')+u_1(s',z').$$
Let us rewrite the last equation with $z'=z$, as that is what $z'$ is. Then
$${\partial g(s',z)\over \partial z}=h(s'-z)-zh'(s'-z)+u_1(s',z).\eqno(4)$$
If 
$${\partial g(s,z)\over \partial z}={\partial g(s',z)\over \partial z}.$$
then, inserting $z'=z$ and $s'=s_0$ gives 
$$h(s)+u(s,z)=h(s_0-z)-zh'(s_0-z)+u_1(s_0,z).$$
The right side is finite when $z=0$ for any pole $s_0$. 
The functions $u(s,z)$ and $u_1(s,z)$ are finite. 
We get a proof of the Riemann Hypothesis: the left side must be 
finite, so $h(s_0)$ is finite. 

But here is an error caused by using the same notation for $z$ and $z'$. 
They are the same number, but they appear in two different roles.
Equation (4) should be calculated as
$${\partial g(s',z')\over \partial z}={\partial z'h(s'-z')\over \partial z}$$
$$={ds'\over dz}{\partial z'h(s'-z))\over \partial s'}+{dz'\over dz}{\partial z'h(s'-z'))\over \partial z'}$$
$$=1\cdot z'h'(s'-z')+1\cdot (h(s'-z')-z'h'(s'-z'))+u_1(s',z')$$
$$=h(s'-z')+u_1(s',z')=h(s'-z)+u_1(s',z)=h(s)+u(s,z).$$

This is what happens in (3): the left side is not the finite partial 
derivative of the function $\varphi_2(s_0,z)$. It is a (potentially) infinite
partial derivative in the coordinates $(s,z)$.  
The confusion was caused by $z'$ being identical to $z$. 

The error is fatal. It was not so obvious in the paper as here. 
I gave several valid proofs of a lemma (and the last proof that was incorrect) 
showing that $\varphi_1(s,z)$ is finite, as it is. I never considered 
the possibility that the left side might not be the partial derivative 
of this finite function, as it does look like it is. 

There are two lessons that can be learned from this false proof. One is
that by setting partial derivatives in a suitable way one can construct
an apparent paradox that may very possibly fool most students, and even
experienced mathematicians. 

The second lesson is more important. In the Spring 2020 I submitted one 
version of [1], with three correct
lemmas and a correct theorem, but with the error in connection of lemma 2
to lemma 3: lemma 3 expected a partial derivative in variables $(s,z)$, while
lemma 2 proved finiteness in variables $(s',z)$. The paper was for 16 weeks
in the Annals of Mathematics. I asked the journal if the paper is in review, 
and received the answer that it is being reviewed. Finally the paper was 
rejected without any referee statement. It could be that the journal 
considered the level of the paper too low and not worthly of a statement, 
but why then keep it for 16 weeks and why to claim that it was in
review? It is possible that the referee of the journal did not spot the error.

This raises the question how to get proof attempts to famous open problems
checked. Submitting them to good, but not leading, journals will almost
certainly stop to the editor, who does not want to deal with such papers.
Submitting them to leading journals usually also stops to the editor, but
even if the paper passes to a referee, there is no guarantee that any referee
statement be given and if given that the statement points out the error. 
It is also possible to submit a paper to journals that are not considered 
good, but if a proof of a famous open problem is published in such a journal, 
it is simply ignored. Thus, this option does not give anything.  
Yet, it is almost certain 
that famous open problems are not solved unless some people make efforts to 
solve them. Possibly a payable review is the only solution, though it is
quite expensive. For instance, I have tried all seven CMI problems. If
checking the solutions requires some 14,000 dollars, that is some money.
Maybe the normal review process could be somewhat improved. The lesson
I get from the 2008 proof attempt is that the review process is not yet
quite ideal. The other lesson, of partial derivatives, well, based on it
you can make a nice trick to try on your students. The new proof of the 
Riemann Hypothesis is in [3].

\end{document}